\documentclass[11pt]{article}
\usepackage{amsfonts,euscript,theorem}
\newcommand{\no}{\noindent}
\newcommand{\bc}{\begin{center}}
\newcommand{\ec}{\end{center}}
\newcommand{\be}{\begin{equation}}
\newcommand{\ee}{\end{equation}}
\newcommand{\bea}{\begin{eqnarray*}}
\newcommand{\eea}{\end{eqnarray*}}
\newcommand{\bean}{\begin{eqnarray}}
\newcommand{\eean}{\end{eqnarray}}
\newcommand{\ran}{\ensuremath{\rangle}}
\newcommand{\lan}{\ensuremath{\langle}}
\newcommand{\rar}{\ensuremath{\rightarrow}}

\newcommand{\pd}{\ensuremath{\partial}}
\newcommand{\pdp}{\ensuremath{\partial}^p}
\newcommand{\pdx}{\ensuremath{\partial}_x}

\newcommand{\h}[1]{\ensuremath{\hbox{ #1 }}}

\newcommand{\w}[1]{\ensuremath{\widetilde{#1}}}

\newcommand{\de}{\ensuremath{\delta}}
\newcommand{\te}{\ensuremath{\theta}}

\newcommand{\Om}{\ensuremath{\Omega}}
\newcommand{\si}{\ensuremath{\sigma}}

\newcommand{\ps}{\ensuremath{\psi}}
\newcommand{\ve}{\ensuremath{\varepsilon}}
\newcommand{\vp}{\ensuremath{\varphi}}
\newcommand{\pe}{\ensuremath{\vp_{\ve}}}
\newcommand{\ka}{\ensuremath{\kappa}}

\newcommand{\Z}{\ensuremath{\mathbb{Z}}}
\newcommand{\N}{\ensuremath{\mathbb{N}}}
\newcommand{\NN}{\ensuremath{\mathbb{N}_0}}
\newcommand{\R}{\ensuremath{\mathbb{R}}}

\newcommand{\C}{\ensuremath{\mathbb{C}}}
\newcommand{\A}{\ensuremath{A_0(\R)}}
\newcommand{\Aq}{\ensuremath{A_q(\R)}}
\newcommand{\G}{\ensuremath{{\EuScript{G}}}}
\newcommand{\GR}{\ensuremath{{\EuScript{G}}(\R)}}

\newcommand{\D}{\ensuremath{{\EuScript{D}}(\R)}}

\newcommand{\DD}{\ensuremath{{\EuScript{D}}'(\R)}}

\theorembodyfont{\itshape}
\newtheorem{Th1}{Theorem}
\newtheorem{Cor1}{ Corollary}
\newtheorem{Cor2}[Cor1]{ Corollary}
\newtheorem{Th2}[Th1]{Theorem}

\begin{document}
\setlength{\baselineskip}{18pt}

\bc {\bf { ON MODELLING OF SINGULARITIES AND THEIR PRODUCTS IN
COLOMBEAU ALGEBRA \boldmath $\GR$ \unboldmath}}

\vspace*{2mm} {\textsc{Blagovest Damyanov}}\footnote{ E-mail:
bdamyanov@mail.bg}

\textit{Bulgarian Academy of Science, INRNE - Theor. Math. Physics
Dept. \\ 72 Tzarigradsko shosse, 1784 Sofia, Bulgaria} \ec

\vspace*{2mm} \setlength{\baselineskip}{13pt} {\small Modelling of
singularities given by discontinuous functions or distributions by
means of generalized functions has proved useful in many problems
posed by physical phenomena.  We introduce in a systematic way
generalized functions of Colombeau that model singularities given
by distributions with singular point support. Moreover, we
evaluate various products of such generalized models whenever the
results admit associated distributions. The results obtained
follow the idea of a well-known result of Jan Mikusi\'nski on
balancing of singular distributional products.  }

\setlength{\baselineskip}{18pt}


\vspace*{2mm}

\no \textbf{1 \ INTRODUCTION}

\vspace*{2mm} The Colombeau algebra of generalized functions \G \
\cite{col84} have become a powerful tool for treating differential
equations with singular coefficients and data \,as well as
singular products of Schwartz distributions. The flexibility of
Colombeau theory allows to model such singularities by means of
appropriately chosen generalized functions, treat them in this
framework and obtain results on distributional level, using the
association process in \G.

In particular, Colombeau functions have proved useful in studying
Euler-Lagrange equations for classical particle in \de-type
potential as well as the geodesic equation for impulsive
gravitational waves; see \cite[\S 1.5, \,\S 5.3]{gct}. Generalized
models in \G \ of Heaviside step-function \te \ were successfully
applied to solving problems arising in Mathematical Physics
\cite{col87}. Other examples involving \te- \,and \de-type
singularities that describe realistic physical phenomena are jump
conditions in hyperbolic systems leading to travelling \de-waves
solutions \cite{col89}, controlled hybrid systems \cite{HO},
geodesics for impulsive gravitational waves  \cite{St}. A detailed
presentation of results on this topic and list of citations can be
found in \cite{mob} and \cite{gct}.

Recall further the well-known result published by Jan Mikusi\'nski
in \cite{mik}\,: \be x^{-1}\,.\,x^{-1} \,-\,
\pi^{{}2}\,\de(x)\,.\, \de(x) \ = \ x^{-2}, \ \ x \in \R.
\label{mik} \ee Though, neither of the products on the left-hand
side here exists, their difference still has a correct meaning in
the distribution space \DD. \ Formulas including balanced singular
products of distributions can be found in the mathematical and
physical literature. For balanced products of this kind, we used
the name `products of Mikusi\'nski type' in a previous paper
\cite{ijpm}, where we derived a generalization in Colombeau
algebra of equation (\ref{mik}) so that the distributions $x^{-
p}$ and $\de^{(q)}$ for arbitrary natural $p, q$ were involved.

Motivated by the afore mentioned works on generalized models in
the algebra \G, we have introduced in a \textit{unified way}
generalized functions of Colombeau that model singularities of
certain type and have additional properties \cite{dam}. The
singularities we considered were given by distributions with
singular support (the complement to the maximal open set where the
distribution is a $C^{\infty}$-function) in a point $x$ on the
real line \R. For $x=0$, such are Dirac \de-function and its
derivatives, Heaviside step function, the non-differentiable
functions  $x_{\pm}^p$, and the distributions $x_{\pm}^{a}, \,a
\in \R\backslash \Z$.

In this paper, we study generalized models in \G \ of the
distributions $x_{\pm}^{- p}, \,p\in\N$ \, and evaluate various
products of such models whenever the result admits associated
distribution. When computed for the canonical embedding of the
distributions in \G, none of the computed here singular products
admits associated distribution.

\vspace*{3mm}
\no \textbf{2 \ NOTATION AND DEFINITIONS}
\vspace*{2mm}

\no \textbf{2.1.} We recall first the basic definitions of
Colombeau algebra $\G(\R)$\,\cite{col84}.

\vspace*{1mm} \textit{Notation 1.} Let \N \ denote the natural
numbers, $\NN = \N \cup \{0\}$, and \,$\de _{ij} = \{\,1$ \,if
\,$i=j, \, = 0$ \,if $i\ne j$\,\}, \,for \,$i, j\in\NN$. Then we
put for arbitrary $q \in \NN$\,:
\[\Aq = \{ \vp(x) \in \D: \int_{\R} x^{j}\,\vp (x)\,dx  = \de _{0j}, \ j = 0, 1,...,q \},\]
where $\D$ is the space of infinitely differentiable functions
with compact support. For $\vp \in\A$ and $\ve > 0$, we will use
the following notation throughout the paper: \ $\pe = \ve
^{-1}\vp(\ve ^{-1}x)$ \ and \ $s \equiv s(\vp) := sup\ \{|x|: \vp
(x)\ne 0 )\}$. Then clearly \ $s(\pe) =\ve s(\vp)$, and denoting
$\si \equiv \si(\vp, \ve) := s(\pe) > 0$, we have $\si := \ve s =
O(\ve)$, as $\ve \rightarrow 0$, for each $\vp \in\A$. \ Finally,
the shorthand notation $\pdx = d/dx$ \,will be used in the
one-dimensional case too.

\vspace*{1mm} \textsc{Definition 1.} Let $\EuScript{E}\,[\R]$ be
the algebra of functions $F(\vp , x): \A \times \R  \rar \C$ that
are infinitely differentiable for  fixed `parameter' \vp. Then the
generalized functions of Colombeau are elements of the quotient
algebra \ \ $\G \equiv\GR = \EuScript{E}_{\mathrm{M}} [\R]\,/ \
\EuScript{I}\,[\R]$. \ Here $\EuScript{E}_{\mathrm{M}}[\R]$ \,is
the subalgebra of `moderate' functions such that for each compact
subset $K$ of \R \,and $p \in\N_0$ \,there is a $q\in \N$ \,such
that, for each $\vp \in \Aq$, \  $\sup_{x \in K}\,|\pdp \,F(\pe,
x)\,| = O(\ve^{-q}), \hbox{ as}  \ \ve \rar 0_+$, where $\pdp$
denotes the derivative of order $p$. \,The ideal
$\EuScript{I}\,[\R]$  of $\EuScript{E}_{\mathrm{M}}[\R]$ consists
of all functions such that for each compact $K \subset\R $ \ and
any $p\in \N_0$ \,there is a $q\in\N$ \ such that, for every $r
\geq q$ and $\vp \in A_r(\R)$, \ $\sup_{x \in K}\,| \pdp \,F(\pe,
x)\,| = O(\ve^{r-q}), \hbox{ as} \ \ve \rar 0_+$.

The differential algebra \GR \ contains the distributions on \R, canonically embedded as a \C-vector subspace \,by the map  \vspace*{-1mm}
\[ i : \DD \rar \,\G : u \mapsto \w{u} = \{\,\w{u}(\vp, x) := (u * \check{\vp})(x)| \,\vp \in\Aq \,\}, \  \h{where} \ \check{\vp}(x) = \vp(-x).  \]

The equality of generalized functions in  \G \ is very strict and a weaker form of equality in the sense of {\em{ association}} is introduced, which plays a fundamental role in Colombeau theory.

\vspace*{1mm}
\no \textsc{Definition} 2 (a) Two generalized functions $F, G \in\GR$ \,are said to be `associated', denoted $F \approx G$, if for some representatives $F(\vp _\ve, x), G(\vp _\ve, x)$ and arbitrary $\ps (x)\in\D$ \,there is a $q\in\NN$, such that for any $\vp (x)\in\Aq, \ \lim_{{}\ve \rar 0_+} \int_{\R} [F(\vp _\ve , x) - G(\vp _\ve, x)] \ps (x)\,dx = 0.$

(b) A generalized function $F\in\GR$ \ is said to be `associated' with a distribution $ u\in \DD$, denoted $F \approx u$, \,if for some representative $F(\vp _\ve, x)$, \,and arbitrary $\ps (x)\in\D$ \,there is a $q\in\NN$, such that for any $\vp (x)\in\Aq, \ \lim_{{}\ve \rar 0_+} \int_{\R} f(\vp _\ve , x) \ps (x)\,dx = \lan u, \ps\ran.$

\vspace*{1mm}
These definitions are independent of the representatives chosen, and the association is a faithful generalization of the equality of distributions. The following relations hold in \G\,:
\be F \approx u \quad  \& \quad F_1 \approx u_1 \  \Longrightarrow \ F + F_1 \approx u + u_1, \quad \pd F \approx \pd u. \label{lin} \ee

\textbf{Remark.} The equation \,$F\approx u$ \,is asymmetric in the sense that the terms cannot be moved over the $\approx$-sign: on the r.h.s. of it there stands a distribution. Of course, its equivalent relation  $F\approx \w{u}$ \ in \G \ \,is symmetric (and can be written as $F - \w{u}\approx 0$ as well). We prefer however the first, simpler and suggesting, notation for the associated distribution.

\vspace*{2mm}
\textbf{2.2.} We next recall the definition of some distributions to be used in the sequel.

\vspace*{1mm}
\textit{Notation 2.} If $a\in\C$ and Re~$a > -1$, denote as usual the locally-integrable functions\,:
\[ x_+^{\,a} =  \{ x^{\,a} \ \h{if} x > 0 , \quad  = 0 \ \h{if} x<0\}, \qquad x_-^{\,a} =  \{  (-x)^{\,a} \ \h{if} x < 0,  \quad = 0 \ \h{if} x >0 \}.  \]
\[ \ln x_+ =  \{ \ln x \ \h{if} x > 0 , \quad  = 0 \ \h{if} x<0\}, \quad \ln x_- =  \{ \ln (-x) \ \h{if} x < 0,  \quad = 0 \ \h{if} x >0 \}.  \]
\[  \ln |x| = \ln x_+ + \ln x_-, \qquad \mathrm{\ln} |x|\,\mathrm{sgn}\,x = \ln x_+ - \ln x_-.\]

\vspace*{1mm}
The distributions $x_{\pm}^{\,a}$ \,are defined for any $a\in\Om := \{a \in \R: \, a \ne -1, -2,\ldots\}$, by setting
\[ x_{+}^{\,a} =  \pd^r \,x_{+}^{\,a +r }(x), \qquad  x_{-}^{\,a} = (-1)^r \,\pd^r \,x_{-}^{\,a +r }(x), \]
where $r\in \NN$  \,is such that $a+r > -1$ \,and the derivatives are in distributional sense.

This definition can be extended also for negative integer values of $a$ by a procedure essentially due to M. Riesz (see \cite[\S \,3.2]{hor}). For each $\ps(x) \in \D, \ a\mapsto\lan x_{+}^a, \ps\ran$ is an analytic function of $a$ on the set $\Om$. The excluded points are simple poles of this function. For any $p\in\NN$, the residue at $a = -p-1$ \,is \ $\lim_{a\rightarrow -p-1}(a+p+1) \ \lan x_{+}^a, \ps\ran = \psi^{(p)}(0)/p!$. Subtracting the singular part, one gets for any $p\in\NN$\,:
\[\lim_{a\rightarrow -p-1} \ \lan x_{+}^a, \ \ps\ran - \frac{1}{p!}\ \frac{\psi^{(p)}(0)}{a+p+1} = - \frac{1}{p!}\int_0^\infty \ln x \,\ps^{(p+1)}\,dx + \frac{\psi^{(p)}(0)}{p!} \sum_{k=1}^{p} \frac{1}{\,k}.\]
The right-hand side of this equation, which is the principal part of the Laurent expansion, was proposed by H\"{o}rmander in \cite{hor} to define the distribution $x_+^{-p-1}$, acting here on the test-function $\ps(x)$. In view of the notation in 2.2, this is equivalent to the following definition of $x_{+}^{-p-1}$ \,for arbitrary $p\in\NN \ (x\in\R)$\,:
\be x_{+}^{-p-1} =  \frac{(-1)^p}{p!} \ \pdx^{p+1} \ln x_+ \ + \  \frac{(-1)^p\,\ka_p}{p!} \ \de^{(p)}(x).\label{x+p}\ee
We have introduced here the shorthand notation \ $ \ka_p :=  \sum_{k=1}^p 1/k  \ (p\in\NN)$; note that $\,\ka_0 = 0$.
Similar consideration leads to the defining equation
\be x_{-}^{-p-1} =  \frac{-1}{p!} \ \pdx^{p+1} \ln x_- \ + \ \frac{\ka_p}{p!} \ \de^{(p)}(x). \label{x-p}\ee
One checks that the distributions $x_{\pm}^{-p}$ satisfy\,:
\[ \pdx \,x_+^{-p} \ = \ -p \ x_+^{-p-1} \ + \ \frac{(-1)^{p}}{p!} \ \de^{(p)}(x), \quad \pdx \,x_-^{-p} \ = \ p \ x_-^{-p-1} \ - \ \frac{1}{p!} \ \de^{(p)}(x). \]
Moreover, it follows immediately that
\be   x_{+}^{-p}\,|_{x \mapsto -x} = x_{-}^{-p} \qquad \hbox{and also} \qquad x_{+}^{-p} \ + \ (-1)^{p} \,x_{-}^{-p} \ = \  x^{-p} \quad (p\in\N), \label{x+-p}\ee
where $x^{-p}$ \,is defined, as usual, as a distributional derivative of order $p$ of $\ln |x|$.

Alternatively, we define the distribution
\be x^{-p}\,\mathrm{sgn}x \ := \  x_{+}^{-p} \ - \ (-1)^{p} \,x_{-}^{-p}  \quad \quad (p\in\NN). \label{altx-p} \ee
Note that $ x^{-p}\,\mathrm{sgn}x \ \neq x^{-p}$ for arbitrary $p\in\NN$; it also differs from the 'odd' and 'even' compositions $|x|^{-p}\,\mathrm{sgn}x := x_+^{- p} + x_-^{- p} = x^{-p}$ for odd natural $p$ and $|x|^{-p} := x_+^{- p} - x_-^{- p} = x^{-p}$ for even $p$.

Recall finally the definition of the distributions $(x {\pm}i 0)^{-p-1}$ for $p\in \NN$\,:
\be (x {\pm}i 0)^{-p-1} := \lim_{y\rar 0_+}  (x {\pm}i y)^{-p-1} = x^{\,-p-1} \ \mp \,\frac{( - 1)^p \ i \,\pi}{p!} \ \de^{(p)}(x), \quad\ x\in\R.   \label{xip}\ee

\vspace*{3mm} \no \textbf{3 \ MODELLING OF SINGULARITIES IN
COLOMBEAU ALGEBRA}  \boldmath \ \G(\R) \unboldmath \vspace*{2mm}

Consider first generalized functions that model the \de-type
singularity in the sense of association, i.e. being associated
with the \de-function. Since there is an abundant variety of such
functions (together with the canonical imbedding $\w{\de}$ in \G \
of the distribution \de), we can put on the generalized functions
in question an additional requirement. So define, following
\cite[\S 10]{mob}, a generalized function $D \in \G$ with the
properties: \be D \approx \de, \qquad D^{\,2} \approx \de.
\label{DD} \ee

To this aim, we let $\vp \in\A$, $s \equiv s(\vp)$, and $\si = s(\pe) = \ve s$ be as in Notation~1, and $D \in \G$ be the class $[ \vp \mapsto D( s(\vp), x)]$. We specify further that $D(s, x) = f(x) + \lambda_{s}\ g(x)$, where $f, g \in \D$ are real-valued, symmetric, with disjoint support, and satisfying:
\[ \int_{\R}f(x) dx = 1, \quad \int_{\R}g(x) dx = 0, \ \h{and} \ \lambda_{s}^2 = \frac{s - \int f^{\,2}(x) dx}{\int g^{\,2}(x) dx}.\]
It is not difficult to check that, for each $\vp \in \A$, the
representative $D(s, x)$ of the generalized function $D$ satisfies
the conditions: \be D(., x)\in \D, \quad D(., -x ) = D(., x),
\quad \frac{1}{s} \int_{\R} D^{\,2}(s, x)\ dx \ = \ \int_{\R} D(s,
x) dx \ = \ 1, \label{D-DD}\ee for each real positive value of the
parameter $s$. Moreover, the generalized function $D$ so defined
satisfies the association relations (\ref{DD}). To show this,
denote by \be D_{\si}(x) \ := \ \frac{1}{\si} D\!\left(\si,
\frac{x}{\si}\right), \  \h{where} \  \si = s(\pe).
\label{strict}\ee Now, for an arbitrary test-function $\ps\in\D$,
\,evaluate the functional values
\[ I_1(\si) = \lan D_{\si}(x), \ps(x) \ran, \qquad  I_2(\si) = \lan D_{\si}^{\,2} (x), \ps(x) \ran, \]
as $\ve \rar 0_+$, or equivalently, as $\si \rar 0_+$. But in view of (\ref{D-DD}), it is immediate to see that \ $\ \lim_{{}\si \rar 0_+} I_1(\si) = \ \lim_{{}\si \rar 0_+} I_2(\si) = \lan \de, \ps \ran$; which according to Definition~2\,(b) gives (\ref{DD}).

The first equation in (\ref{DD}) is in consistency with the observation that $ D_{\si}(x) $ \,is a strict \de-net \,as defined in distribution theory \cite[\S 7]{mob}. But notice that $D$ is not the canonical embedding $\w{\de}$ of the \de-function since ${\w{\de}}^{\,2}$ does not admit associated distribution.

The flexible approach to modelling singularities allowed by generalized functions in \G, so that the models satisfy auxiliary conditions, can be systematically applied to defining generalized models of particular singularities. We will consider models of singularities given by distributions with singular point support. For their definition, we intend to take advantage of the properties of \de-modelling function $D$. Observe that it holds
\[ (\de*D(s, .))(x)  = \lan \de_y, \,D(s, x-y) \ran  =  D(s, x). \]
\[ \vspace*{-1mm} (\de\,'* D(s, .))(x)  = \lan \de\,'_y, \,D(s, x-y) \ran   = - \lan \de_y\, \,\pd_y D(s, x-y) \ran  = \lan \de_y, \,D\,'(s, x-y) \ran  =  D\,'(s, x). \]

\no This can be continued by induction for any derivative to define a generalized function $D^{\,(p)}(x)$ \,that models $\de^{\,(p)}(x)$ and has representative \ $D^{\,(p)}(s, x) = (\de^{\,(p)}* D(s, .))(x)$.

\no Clearly, this definition is in consistency with the differentiation: \ $\pdx D^{\,(p)}(x) = D^{\,(p+1)}(x)$, $p\in\NN$. Moreover,
\be D^{\,(p)}(-x) \ = \ (-1)^p D^{\,(p)}(x). \label{D-}\ee

In \cite{dam} we have employed such procedure for a unified modelling of singularities  given by distributions with singular point support, i.e. (besides $\de^{(p)}$) the distributions $ x_{\pm}^a, \,a\in\Om$. Namely, choosing an arbitrary generalized function $D$ with representative $D(s,x)$ that satisfies (\ref{D-DD}) for each $\vp \in \A$, we have introduced generalized functions $ X_{\pm}^{a}(x) $, modelling the above singularities, with representatives
\be X_{\pm}^{a}(s, x) :=  ( y_{\pm}^{a}*D(s, y))(x), \quad a\in\Om.\label{defXa}\ee
This is consistent with the differentiation\,: \,$\pdx X_{\pm}^{a}(x)  =  a X_{\pm}^{a-1}(x)$; \,in particular, $H'  = D$, \,where $H\in\G$ is model of the step-function $\te$, with representative $H(s,x) = \te *D(s, .)(x)$.

Extending now definition (\ref{defXa}) to the distributions $ x_{\pm}^{- p-1}, \,p\in\NN$, we can write
\be X_{\pm}^{- p -1}(s, x) :=  ( y_{\pm}^{- p-1}*D(s, y))(x). \label{Xpm-}\ee
Similarly, we put \ $ \mathrm{Ln}\,x_{\pm} := \ln y_{\pm} * D(s, y))(x)$.

We note that generalized functions so introduced are indeed models of the corresponding singularities. It is not difficult to show that for each $ a\in \Om $
\[ X_{\pm}^{a}(x) \ \approx \ x_{\pm}^{a}(x); \ \h{in particular,} \,H  \approx \te, \ \hbox{and} \ \,H^p \approx \te \ \hbox{for each} \ p\in\N.\]
Also, it was proved in \cite{dam} that --- as it can be expected --- the functions $H$ and $D$ that model correspondingly the \te- and \de-type singularities satisfy the relation \ $H\,.\,D \ \approx \ \frac{1}{2}\ \de$. \ Moreover, these generalized models were proved to satisfy
\be H\,.\,D\,' \ \approx \ - \ \de \ + \ \frac{1}{2}\ \de\,'. \label{HD'}\ee
Concerning the singularities given by the distributions $x_{\pm}^{\,-p}, \ p\in\N$, it can be easily checked that $ \mathrm{Ln}_{\pm}x \approx \ln_{\pm}x $  for the latter locally-integrable function. Then the modelling property for the generalized functions $ X_{\pm}^{- p }(x)$ \,follows in view of relation (\ref{lin}) for consistency between differentiation and association in \G.

\vspace*{1mm}
Finally, we need to compute the representatives of the generalized models when they depend on $\pe$, or rather on the value $s(\pe) = \ve  s(\vp) = \si$. In view of equations (\ref{x+p}),  (\ref{x-p}), (\ref{strict}), (\ref{defXa}), and (\ref{Xpm-}), we obtain for the corresponding representatives $(p\in \NN)$\,:
\be X_{+ \,\si}^{p}(x) = \frac{1}{\si} \int_0^{\infty} y^p D\!\left(\si, \frac{x-y}{\si}\right)\,dy, \quad X_{- \,\si}^{p}(x) = \frac{1}{\si} \int_{- \infty}^0 (- y)^p D\!\left(\si, \frac{x-y}{\si}\right)\,dy.  \label{Xpsi}\ee
\[ X_{+ \,\si}^{- p-1}(x) = \frac{(- 1)^p}{\si^{p+2}\ p!} \int_0^{\infty} \ln y \,D^{(p+1)}\left(\si, \frac{x-y}{\si}\right)\,dy + \,\frac{(- 1)^p\,\ka_p}{\si^{p+1}\,p!}\,D^{(p)}\left(\si, \frac{x}{\si}\right).  \]
\be X_{- \,\si}^{- p-1}(x) = \frac{ - 1}{\si^{p+2}\ p!} \int_{- \infty}^0 \ln (- y) \,D^{(p+1)}\left(\si, \frac{x-y}{\si}\right)\,dy + \,\frac{\ka_p}{\si^{p+1}\,p!}\,D^{(p)}\left(\si, \frac{x}{\si}\right). \label{X-psi}\ee

\vspace*{5mm} \no \textbf{4 \ PRODUCTS OF SOME SINGULARITIES
MODELLED IN} \boldmath \ \G(\R) \unboldmath

\vspace*{2mm}
The models of singularities we consider all have products in Colombeau algebra as generalized functions, but we are seeking results that can be evaluated back in terms of distributions, i.e. such that admit an associated distribution. We will establish first certain balanced products of generalized models in the algebra \GR\ that exist on distributional level, proving the following.

\no \begin{Th1}$\!\!.$   The generalized models of the distributions $x_{\pm}^{- 2},  \te, \check{\te}$, and $\de\,'(x)$ satisfy:
\be
X_-^{- 2}\,.\,H -  \mathrm{Ln}\,x_+\ .\,D\,' \ \approx \ - \ \de. \label{th1-}\ee
\be X_+^{- 2}\,.\,\check{H}  + \mathrm{Ln}\,x_- \ .\,D\,' \ \approx \ - \ \de. \label{th1+}
\ee
\end{Th1}

\no \textit{Proof:} (i) For an arbitrary test-function $\ps(x) \in \D$, denote \ $I(\si) := \lan\  X_{- \,\si}^{- 2}\,.\,H_{\si}, \,\ps(x)\,\ran $. \,Suppose (without loss of generality) that \ supp~$\!D(\si,x) \subseteq [-l, l]$ \,for some $l\in \R_+$; then $-l\leq x/\si \leq l$ implies \ $-l\si \leq x \leq l\si$. \,Now from equations (\ref{X-psi}) for $p=1$ and (\ref{Xpsi}) for $p=0$,  we get on transforming the variables  $y=  \si u + x, \,z =  \si v + x $, and \,$x = - \si w$
\bean I(\si) & = & - \,\frac{1}{\si^4} \int_{-\si l}^{\,\si l} dx\,\ps(x) \int_0^{\si l +x} dy\,D\!\left(\si, \frac{x-y}{\si}\right) \int_{- \si l + x}^0 \ln (- z) \,D\,'' \left(\si, \frac{x-z}{\si}\right)\,dz \nonumber \\ & &
+ \,\frac{1}{\si^3} \int_{-\si l}^{\,\si l} dx\,\ps(x) D\,' \left(\si, \frac{x}{\si}\right) \int_0^{\si l +x}D\!\left(\si, \frac{x-y}{\si}\right)\,dy \nonumber \\ & = &
- \,\frac{1}{\si}\int_{- l}^{\,l} dw\,\ps(- \si w) \int_w^{\,l} du\,D(\si, u) \int_{- l}^w \ln (\si w - \si v) \,D\,''(\si, v)\,dv \nonumber \\ & &
+ \,\frac{1}{\si} \int_{- l}^{\,l} dw\,\ps(- \si w) D\,' (\si, w) \int_w^{\,l}D(\si, u)\,du  \ =: \ I_1 + I_2. \label{I12}\eean
Applying Taylor theorem to the test-function \ps \ and changing the order of integration (which is permissible here), we get further
\bea I_1 & = & - \,\frac{\ps(0)}{\si}\int_{- l}^{\,l} du\,D(\si, u) \int_{- l}^u dv\,D\,''(\si, v) \int_v^{\,u} \ln (\si w - \si v) \,dw \\ & &
+ \,\ps\,'(0) \int_{- l}^{\,l} du\,D(\si, u) \int_{- l}^u dv\,D\,''(\si, v) \int_v^{\,u} \ln (\si w - \si v) \,w \ dw \ + \ o\,(1) \eea
Here the Landau symbol $ o(1)$ stands for an arbitrary function of asymptotic order less than any constant, and the asymptotic evaluation is obtained taking into account that the third term in the Taylor expansion is multiplied by definite integrals majorizable by constants. \,Now the substitution $w \,\rar \,t = (w-v)/(u-v)$, together with $w-v = (u-v) t$, yields
\bea I_1 & = & - \,\frac{\ps(0)}{\si}\int_{- l}^{\,l} du\,D(\si, u) \int_{- l}^u dv\,D\,''(\si, v) (u-v) \left[\ln (\si u - \si v) + \int_0^1 \ln t \,dt\right] \\ & &
+ \,\ps\,'(0) \int_{- l}^{\,l} du\,D(\si, u) \int_{- l}^u dv\,D\,''(\si, v) (u-v)^{\,2} \left[\frac{1}{2} \ln (\si u - \si v) + \int_0^1 t \,\ln t \,dt\right] \\ & &
- \,\ps\,'(0) \int_{- l}^{\,l} du\,D(\si, u) \int_{- l}^u dv\,D\,''(\si, v) (u-v)^{\,2} \left[\ln (\si u - \si v) + \int_0^1 \ln t \,dt\right] \\ & &
+ \,\ps\,'(0) \int_{- l}^{\,l} du\,u\,D(\si, u) \int_{- l}^u dv\,D\,''(\si, v) (u-v) \left[\ln (\si u - \si v) + \int_0^1 \ln t \,dt\right] \ + \ o\,(1).\eea
Calculating the integrals \,$\int_0^1 \ln t \,dt = -1$, $\int_0^1 \ln t \,dt = -1/4$, replacing $v = u - (u-v)$, and integrating by parts in the variable $v$ (the integrated part being 0) we get
\bean I_1 & = & - \,\frac{\ps(0)}{\si}\int_{- l}^{\,l} du\,D(\si, u) \int_{- l}^u \ln (\si u - \si v) \ D\,'(\si, v) dv  \ - \ 2\,\ps(0)\nonumber \\ & &
+ \,\ps\,'(0) \int_{- l}^{\,l} du\,u\,D(\si, u) \int_{- l}^u \ln (\si u - \si v)\,D\,'(\si, v) \,dv \nonumber \\ & &
- \,\ps\,'(0) \int_{- l}^{\,l} du\,D(\si, u) \int_{- l}^u \ln (\si u - \si v) \,D(\si, v)\,dv  \ + \ o\,(1).\label{I1}\eean
To obtain the latter result, we have used equation (\ref{D-DD}) and also that
\be \int_{- l}^{\,l} du\,D(\si, u) \int_{- l}^u D(\si, v)\,dv = \frac{1}{2}. \label{prop}\ee

Applying again Taylor theorem to the function \ps, changing the order of integration, and integrating by parts in the variable $w$, we obtain for the second term in (\ref{I12})\,:
\bea  I_2 & =  & \!\!\!\frac{\ps(0)}{\si} \int_{- l}^{\,l} dw\, D\,' (\si, w) \int_w^{\,l}D(\si, u)\,du  -  \ps\,'(0) \int_{- l}^{\,l} dw\,w\ D\,' (\si, w) \int_w^{\,l}D(\si, u)\,du \\ & = &
\!\!\!\frac{\ps(0)}{\si}\!\int_{- l}^{\,l} du\,D(\si, u) \!\int_{- l}^u D\,' (\si, w)\,dw -  \ps\,'(0) \!\int_{- l}^{\,l} du\,D(\si, u) \!\int_{- l}^u \,w\ D\,' (\si, w)\,dw  + o\,(1) \\ & = &
\ps(0) \ - \ \frac{1}{2}\  \ps\,'(0) \  + \ o\,(1), \eea
where equation (\ref{prop}) is used again. \ Therefore
\bean I(\si) \!\!\! & = & \!\!\!- \frac{\ps(0)}{\si}\!\int_{- l}^{\,l}\!du\,D(\si, u) \!\int_{- l}^u \ln (\si u - \si v) D\,'(\si, v)\,dv \nonumber \\ & & \!\!\!+ \ps\,'(0) \!\int_{- l}^{\,l}du\,u\,D(\si, u)\!\int_{- l}^u \ln (\si u -\si v)\,D\,'(\si, v)\,dv \nonumber \\ & &
\!\!\!- \ps\,'(0)\!\int_{- l}^{\,l} du\,D(\si, u)\!\int_{- l}^u \ln (\si u - \si v)\,D(\si, v)\,dv -  \ps(0) -  \frac{1}{2} \ps\,'(0)  + o\,(1).\label{Isi}\eean

(ii) On the other hand, denoting \ $J(\si) := \lan \,\mathrm{Ln}\,x_{+ \si}\,.\,D\,'_{\si}, \,\ps(x)\,\ran$, we obtain on transforming the variables  $y= \si u + x$ \,and \,$x = - \si v$, applying Taylor theorem to the test-function \ps, \,and changing the order of integration\,:
\bea J(\si) & = & \,\frac{1}{\si^3} \int_{-\si l}^{\,\si l} dx\,\ps(x) \,D\,' \left(\si, \frac{x}{\si}\right) \int_0^{\si l +x} \ln y \ D\!\left(\si, \frac{x-y}{\si}\right) dy \\ & = &
- \,\frac{1}{\si}\int_{- l}^{\,l} dv\,\ps(- \si v) \,D\,'(\si, v)\int_v^{\,l} \ln (\si v - \si u) \,D(\si, u)\, du  \nonumber \\ & = &
- \,\frac{\ps(0)}{\si}\int_{- l}^{\,l} du\,D(\si, u) \int_{- l}^u \ln (\si v - \si u) \,D\,'(\si, v)\,dv \\ & &
+ \,\ps\,'(0) \int_{- l}^{\,l} du\,D(\si, u) \int_{- l}^u \ln (\si v - \si u) \,v\,D\,'(\si, v)\,dv \ + \ o\,(1).\eea
Replacing now $v = u + (v - u)$ in the last term and integrating by parts the third term so obtained, we get
\bean J(\si) & = &
- \,\frac{\ps(0)}{\si}\int_{- l}^{\,l} du\,D(\si, u) \int_{- l}^u \ln (\si v - \si u) \,D\,'(\si, v)\,dv \nonumber \\ & &
+ \,\ps\,'(0) \int_{- l}^{\,l} du\,u\,D(\si, u) \int_{- l}^u \ln (\si v - \si u) \,D\,'(\si, v)\,dv \nonumber  \\ & &
-  \,\ps\,'(0) \int_{- l}^{\,l} du\,D(\si, u) \int_{- l}^u \ln (\si v - \si u)\,D(\si, v)\,dv - \frac{1}{2}\,\ps\,'(0) \ + \ o\,(1). \label{Jsi}\eean

Combining now equations (\ref{Isi}) and (\ref{Jsi}), we obtain by linearity
\[\lim_{{}\si \rar 0_+} \int_{\R} \ps(x) \left[\,X_{- \,\si}^{\,- 2}(x)\,.\,H_{\si}(x) - \,\mathrm{Ln}\,x_{+ \,\si}(x)\,.\,D\,'_{\si}(x) \right] \,dx = - \,\ps(0)  = - \,\lan \de, \ps\ran.
\]
According to Definition~2(b), this proves the first equation in (\ref{th1-}). The second equation follows on replacing $x \rightarrow - x$ in the first one and taking into account equations (\ref{x+-p}) and (\ref{D-}). This completes the proof.

\vspace*{2mm}
The above balanced products of the functions $X_{\pm}^{\,- 2}$ supported in the corresponding real half-lines can be employed further to get results on singular products of the generalized modelling functions $ X^{\,- 2}\,\mathrm{sgn}\,x $ and $X^{\,- 2}$ (obtained from equations (\ref{altx-p}), (\ref{x+-p}), and (\ref{Xpm-})).

\begin{Cor1}$\!\!.$ The following balanced product holds for the generalized models of the distribution $x^{\,- 2}\,\mathrm{sgn}\,x $, \te, and $\de\,'$\,:
\be X^{\,- 2}\,{\mathrm{sgn}}\,x \,.\,H  \ + \ {\mathrm{Ln}}\,|x|\,{\mathrm{sgn}}\,x \,.\,D\,' \ \approx \ x_+^{\,-2} + 2\,\de. \label{cor1H}\ee
\end{Cor1}

\no \textit{Proof:} Consider the following chain of identities and associations in \GR, taking into account equation (\ref{th1+}) and the relation \ $H + \check{H} \approx 1$\,:
\[ X_+^{\,- 2}\, .\,H \,= \, X_+^{\,- 2}\, .\,( 1 - \check{H}) \, = \, X_+^{\,- 2}\, - \, X_+^{\,- 2}\,.\,\check{H}  \ \approx \ X_+^{\,- 2}\, + \ \mathrm{Ln}\,x_- \ .\,D\,' \ + \ \de.\]
Thus
\[ X_+^{\,- 2}\, .\,H \, - \ \mathrm{Ln}\,x_- \ .\,D\,' \ \approx  \ \ X_+^{\,- 2} \ + \ \de,\]
which, in view of the association \ $ X_+^{\,- 2} \approx x_+^{\,- 2}$ and the linearity by (\ref{lin}) of the association in \G \, leads to the balanced product
\be X_+^{\,- 2}\, .\,H \, - \ \mathrm{Ln}\,x_- \ .\,D\,' \ \approx x_+^{\,-2} \ + \ \de. \label{cor1+} \ee

Further, equation (\ref{altx-p}) for $p=2$, as well as equations (\ref{th1-}) and (\ref{cor1+}), will all yield
\[ X^{\,- 2}\,{\mathrm{sgn}}\,x \, .\,H \, \,= \,\left(X_+^{\,- 2} \ - \ X_-^{\,- 2}\right)\, .\,H \, \,\approx \ \mathrm{Ln}\,x_- \ .\,D\,' \ + \ x_+^{\,-2} +  \de - \mathrm{Ln}\,x_+\ .\,D\,'  + \de.\]
In view of relation  (\ref{lin}) for linearity of the association, this proves equation (\ref{cor1H}).

\vspace*{2mm}
Other consequences from the above results are given by this.

\no \begin{Cor2}$\!\!.$  The generalized models in \G \ of the distributions $(x {\pm}i 0)^{\,- 2}$, $\te$, and $\de\,'$ \ satisfy
\be ( X {\pm}i 0)^{\,- 2}\,.\,H \ - \ {\mathrm{Ln}}\,|x|\ .\,D\,' \ \approx  \ x_+^{\,- 2} \  \mp \ i \pi\ \de(x) \pm \ \frac{ i \pi }{2}\,\de\,'. \label{cor2}\ee
\end{Cor2}

\textbf{Proof}\,:  The second equation in (\ref{x+-p}), as well as equations (\ref{th1-}) and (\ref{cor1+}), now give
\[ X^{\,- 2} \, .\,H \, \,= \,\left(X_+^{\,- 2} \ + \ X_-^{\,- 2}\right)\, .\,H \, \,\approx \ \mathrm{Ln}\,x_- \ .\,D\,' \ + \ x_+^{\,-2} +  \de + \mathrm{Ln}\,x_+\ .\,D\,'  -  \de.\]
In view of (\ref{lin}), this yields
\be X^{\,- 2} \,.\,H  \ - \ \mathrm{Ln}\,|x| \,.\,D\,' \ \approx \ x_+^{\,-2}. \label{cor2+}\ee

Employing further equations (\ref{xip}), (\ref{cor2+}) and (\ref{HD'}), we get
\[ ( X {\pm}i 0)^{-1}\, .\,H \, = \,X^{\,- 2}\, .\,H \ \pm \ i \pi \,D\,'(x)\,.\,H \ \approx \mathrm{Ln}\,|x| \,.\,D\,' + \ x_+^{\,- 2} \ \mp \ i \pi\ \de \pm \ \frac{ i \pi }{2}\,\de\,',\]
which in view of linearity of association in \G \ proves equation (\ref{cor2}).

\vspace*{2mm}
Finally, we will evaluate some products of singularities given by the non-differentiable functions $x_{\pm} $ modelled by the generalized functions $X_{\pm} $ with derivatives of $D$. They only exist as balanced products, as demonstrated by this.

\no \begin{Th2}$\!\!.$   The following balanced products hold for the modelling generalized function $X_{\pm}, H $ and $D$ \,:
\be X_+\,.\,D^{\,(4)} \ + \ H\,.\,D^{\,(3)} \ \approx \ \frac{5}{2}\,\de\,'' \ - \ \frac{3}{2}\,\de\,''' \label{th2+}\ee
\be X_-\,.\,D^{\,(4)} \ + \ \check{H}\,.\,D^{\,(3)} \ \approx \ \frac{5}{2}\,\de\,'' \ + \ \frac{3}{2}\,\de\,''' \label{th2-}\ee
\end{Th2}

\no \textit{Proof:} For an arbitrary $\ps(x) \in \D$, we denote \ $I(\si) := \lan\,X_{+ \ \si}(x)\,.\,D_{\si}^{\,(4)}(x), \,\ps(x)\,\ran $. From equations (\ref{strict}) and (\ref{Xpsi}),  we get on transforming the variables  $y= \si v + x, \,x= - \si u$, \,changing the order of integration, and applying Taylor theorem
\bea
I(\si)  & = & \frac{1}{\si^{\,3}} \int_{-l}^{\,l} du \,\ps(- \si u) D^{\,(4)}(\si, u) \int_u^{\,l} (v-u) D (\si, v)\,dv  \\
& = & \frac{\ps(0)}{\si^{\,3}} \int_{-l}^{\,l} dv\,D(\si, v) \int_{- l}^{\,v} (v-u) D^{\,(4)}(\si, u)\,du  \\
&  & - \,\frac{\ps\,'(0)}{\si ^{\,2}}\ \int_{-l}^{\,l} dv\,D(\si, v) \int_{-l}^{\,v} u\,(v-u) D^{\,(4)} (\si, u)\,du \\
&  & + \,\frac{\ps\,''(0)}{2 \,\si }\ \int_{-l}^{\,l} dv\,D(\si, v) \int_{-l}^{\,v} u^{\,2}\,(v-u) D^{\,(4)}(\si, u)\,du \\ & &
 -  \,\frac{\ps\,'''(0)}{6}\ \int_{-l}^{\,l} dv\,D(\si, v) \int_{-l}^{\,v} u^{\,3}\,(v-u) D^{\,(4)}(\si, u)\,du +  O(\si) \\ & =: & \ps(0) \ I_0 \ + \ \ps\,'(0) \ I_1 \ + \ \ps\,''(0) \ I_2 \ + \ \ps\,'''(0) \ I_3 \ + \ O(\si).
\eea

Denote further \,$J(\si) := \lan\,H_{\ \si}(x)\,.\,D_{\si}^{\,(3)}(x), \,\ps(x)\,\ran$. Proceeding as above, we get
\bean
J(\si)  & = & - \,\frac{\ps(0)}{\si^{\,3}} \int_{-l}^{\,l} dv\,D(\si, v) \int_{- l}^{\,v}   D^{\,(3)}(\si, u)\,du   \nonumber \\
&  & +  \,\frac{\ps\,'(0)}{\si ^{\,2}}\ \int_{-l}^{\,l} dv\,D(\si, v) \int_{-l}^{\,v} u\,  D^{\,(3)} (\si, u)\,du \nonumber \\
&  & - \,\frac{\ps\,''(0)}{2 \,\si }\ \int_{-l}^{\,l} dv\,D(\si, v) \int_{-l}^{\,v} u^{\,2}\, D^{\,(3)}(\si, u)\,du \nonumber \\ & &
+ \,\frac{\ps\,'''(0)}{6}\ \int_{-l}^{\,l} dv\,D(\si, v) \int_{-l}^{\,v} u^{\,3}\,(v-u) D^{\,(3)}(\si, u)\,du +  O(\si) \nonumber \\ & =: & \ps(0) \ J_0 \ + \ \ps\,'(0) \ J_1 \ + \ \ps\,''(0) \ J_2 \ + \ \ps\,'''(0) \ J_3 \ + \ O(\si).\label{Js}\eean
We next compute the terms $I_k, \,k= (0, 1, 2, 3)$. During the calculations we shall use equations (\ref{D-DD}), (\ref{prop}), as well as that
\[ \frac{1}{\si} \int_{- l}^{\,l} v\,D(\si, v) D'(\si, v)\,dv \ = \ - \,\frac{1}{2 \,\si} \int_{- l}^{\,l} D^{\,2}(\si, v)\,dv \ = \ - \,\frac{1}{2}. \]
Also, due to the equality \,$ D\,'(., -x ) = - D\,'(., x)$, the following equations hold
\[ \int_{- l}^{\,l} D(\si, v) D'(\si, v)\,dv \ = \ \int_{- l}^{\,l}v\ D^{\,2} (\si, v) \,dv \ = \ \int_{- l}^{\,l}v^{\,2}\ D(\si, v) D'(\si, v)\,dv \ = \ 0. \]

\no Integrating now by parts in the variable $u$, the integrated part being 0 each time, we obtain\,:
\bea  I_0 & = & \frac{1}{\si^{\,3}} \int_{-l}^{\,l} dv\,D(\si, v) \int_{- l}^{\,v} D^{\,(3)}(\si, u)\,du \ = \ - \,J_0. \\*[1mm]
I_1 & = & - \,\frac{1}{\si ^{2}}\!\int_{-l}^{\,l} dv\,D(\si, v) \!\int_{-l}^{\,v} u\,D^{(3)}(\si, u)\,du + \!\frac{1}{\si ^{2}}\!\int_{-l}^{\,l}\!dv\,D(\si, v)\int_{-l}^{\,v}\!(v - u) D^{(3)} (\si, u)\,du \\ & = &
- J_1 \ + \ \frac{1}{\si ^{\,2}}\!\int_{-l}^{\,l}\!D(\si, v)D^{\,'}(\si, v)\,dv = \ - J_1. \\*[1mm]
I_2 & = & \frac{1}{2 \si}\!\int_{-l}^{\,l} dv\,D(\si, v)\!\int_{-l}^{\,v} u^{\,2}\,D^{\,(3)}(\si, u)\,du - \frac{1}{\si}\!\int_{-l}^{\,l} dv\,D(\si, v)\!\int_{-l}^{\,v} u (v-u)\,D^{\,(3)}(\si, u)\,du \\ & = &\ - \ J_2  \ + \ I\,'_2, \eea

\vspace*{-3mm}
\no where

\vspace*{-8mm}
\bea I\,'_2 & = & \!\!\!\frac{1}{\si}\!\int_{-l}^{\,l}\!dv\,D(\si, v)\!\int_{-l}^{\,v}\! (v-u)^{2}D^{(3)}(\si, u)\,du - \!\frac{1}{\si}\!\int_{-l}^{\,l}\!dv\,v\,D(\si, v)\!\int_{-l}^{\,v}\! (v-u)D^{(3)}(\si, u)\,du \\ & = &
\!\!\!\frac{2}{\si}\!\int_{-l}^{\,l} D^{\,2}(\si, v)\,dv  - \frac{1}{\si}\!\int_{-l}^{\,l} v \ D(\si, v)\,D\,'(\si, v)\,dv \ = \ \frac{5}{2}.\\*[1mm]
I_3 & = & - \!\!\!\frac{1}{6}\ \int_{-l}^{\,l} dv\,D(\si, v)\!\int_{-l}^{\,v} u^{\,3}\,D^{\,(3)}(\si, u)\,du + \!\frac{1}{2} \!\int_{-l}^{\,l} dv\,D(\si, v)\!\int_{-l}^{\,v} u^{\,2}\,(v-u) D^{\,(3)}(\si, u)\,du \\ & = & - \ J_3 \ + \ I\,'_3. \eea

\vspace*{-3mm}
\no Here

\vspace*{-8mm}
\bea I\,'_3 & = & \frac{1}{2} \!\int_{-l}^{\,l}\!dv\,D(\si, v)\!\int_{-l}^{\,v}\!(v-u)^{3}\,D^{(3)}(\si, u)\,du + \!\int_{-l}^{\,l}\!dv\,v\,D(\si, v)\!\int_{-l}^{\,v}\!u\,(v-u)\,D^{(3)}(\si, u)\,du \\ & &
-  \frac{1}{2} \!\int_{-l}^{\,l} dv\,v^{\,2}\,D(\si, v)\!\int_{-l}^{\,v} (v-u)\,D^{\,(3)}(\si, u)\,du  = 3 \int_{- l}^{\,l} dv\,D(\si, v) \int_{- l}^v D(\si, u)\,du \\ &  &
- \int_{-l}^{\,l}\!dv\,vD(\si, v)\!\int_{-l}^{\,v}\!(v-u)^{2}D^{(3)}(\si, u)\,du +
\!\int_{-l}^{\,l}\!dv\,v^{2}D(\si, v)\!\int_{-l}^{\,v}\!(v-u)D^{(3)}(\si, u)\,du \\ & = & \frac{3}{2}- 2  \int_{- l}^{\,l}v\ D^{\,2} (\si, v) \,dv  + \int_{- l}^{\,l}v^{\,2}\ D(\si, v) D'(\si, v)\,dv \ = \ \frac{3}{2}.\eea
Summing up, we get
\be I(\si) = - \ps(0) \ J_0 \ - \ \ps\,'(0) \ J_1 \ - \ \ps\,'(0) \ J_2 \ - \ \ps\,''(0) \ J_3 \ +  \ \frac{5}{2}\ \ \ps\,''(0) \ + \ \frac{3}{2}\ \ps\,'''(0) \ + \ O(\si). \label{If}\ee
Now from equations (\ref{Js}) and (\ref{If}), we obtain by linearity
\[\lim_{{}\si \rar 0_+} \int_{\R} \ps(x) \left[\,X_{+ \ \si}(x)\,.\,D_{\si}^{\,(4)}(x) + H_{\ \si}(x)\,.\,D_{\si}^{\,(3)}(x)\right] \,dx =  \,\lan  \,\frac{5}{2}\,\de\,''  -  \frac{3}{2}\,\de\,''', \,\ps\,\ran.\]
According to Definition~2(b), this proves equation (\ref{th2+}), whereas equation (\ref{th2-}) follows on replacing $x \rightarrow - x$ in the former. The proof is complete.

\vspace*{2mm}
\textbf{Remark.} When computed for the canonical embedding of the distributions in \G, none of the above singular products can be balanced so as to admit associated distribution.

\vspace*{1mm}
\setlength{\baselineskip}{15pt}

\end{document}